\documentclass{amsart}
\title{A Gray path on binary partitions}
\author{Thomas Colthurst \and Michael Kleber}
\date{sometime between April and July 2003}
\newtheorem{thm}{Theorem}
\newtheorem*{defn}{Definition}
\newcommand{\B}{\mathcal B}
\newcommand{\rdots}{{\stackrel{\leftarrow}{\ldots}\,}}
\newcommand{\ntwo}{\mbox{\small$\frac{n}{2}$}}
\begin{document}
\maketitle

A {\em binary partition} of a positive integer $n$ is a partition of
$n$ in which each part has size a power of two.  Let $b(n)$ denote the
number of binary partitions of $n$.  Since a binary partition either
has a part of size 1 or else is twice a binary partition of $n/2$,
$b(n)$ satisfies the recurrence
\begin{equation}
\label{bn}
\begin{array}{ll}
b(n) = b(n-1),  & n \mbox{ odd;} \\[4pt]
b(n) = b(n-1) + b(n/2), & n \mbox{ even.}
\end{array}
\end{equation}

The generating function $\sum b(n) x^n = 1/\prod (1-x^{2^k})$ and
small values of $b(n)$ were written down by Euler~\cite[\S50]{euler},
and Mahler~\cite{mahler} and de Bruijn~\cite{deBruijn} gave
increasingly good asymptotics for $b(n)$ and for partitions into
powers of $r$ other than $2$; this is sometimes called ``Mahler's
partition problem.''

In this note we first use a variation on the recurrence~(\ref{bn}) to
construct a Gray sequence on the set of binary partitions themselves.
This is an ordering of the set of binary partitions of each $n$ (or of
all $n$) such that adjacent partitions differ by one of a small set of
elementary transformations; here the allowed transformatios are
replacing $2^k+2^k$ by $2^{k+1}$ or vice versa (or addition of a new
$+1$).  Next we give a purely local condition for finding the
successor of any partition in this sequence; the rule is so simple
that successive transitions can be performed in constant time.
Finally we show how to compute directly the bijection between $k$ and
the $k$th term in the sequence.

Thanks to Donald Knuth for requesting\footnote{open problem 
Ex.~59 in TAOCP 7.2.1.4 pre-fascicle 2D, 28 November 2002}
Theorem~\ref{thm_B}, and to Richard Stanley for responding to an
early draft by asking about the possibility of Theorem~\ref{thm_Bk}.

\section{Construction of Gray sequences}

\begin{thm}
\label{thm_B}
For each positive integer $n$, the binary partitions of $n$ can be
arranged in a sequence $\B(n)$ such that adjacent partitions differ by
an operation of the form
$$
\cdots+2^k+2^k+\cdots \,\,\,\longleftrightarrow\,\,\, \cdots+2^{k+1}+\cdots.
$$
Moreover the sequence runs first through all partitions ending with
$\cdots+1$, and next through all partitions ending with $\cdots+2$.
\end{thm}

\proof
We give a recursive construction for $\B(n)$.  We let $Q(n)$ and
$S(n)$ denote the first and last partitions in the sequence (reserving $R$
for later use); $Q(n)$ will always be $1+1+\cdots+1$.
\begin{equation}
\label{Bn}
\renewcommand{\arraystretch}{1.4}
\begin{array}{rl}
n \equiv 1 \bmod 2: & Q(n-1)+1,\ldots,S(n-1)+1 \\
n \equiv 0 \bmod 4: & Q(n-1)+1,\ldots,S(n-1)+1, Q(\ntwo)\times2,\ldots,S(\ntwo)\times2 \\
n \equiv 2 \bmod 4: & Q(n-1)+1,\ldots,S(n-1)+1, S(\ntwo)\times2,\rdots,Q(\ntwo)\times2
\end{array}
\end{equation}
Here $\times2$ indicates doubling each part of a partition.  Note that
in the $2 \bmod 4$ case, the $\ntwo$ sequence appears in reverse
order.

The same logic that justified formula~(\ref{bn}) for $b(n)$ shows that
the above sequences certainly contain all binary partitions.  By induction, 
we only need to check that the two points of concatenation obey the 
adjacency condition.

For $n \equiv 0 \bmod 4$, $S(n-1) +1 = S(n-2)+1+1 = Q(\ntwo - 1)\times2 +1+1$
(since $n-2\equiv 2\bmod 4$), which is $2+\cdots+2+1+1$, since $Q(m)$ is the
all 1's partition for all $m$.  Meanwhile, $Q(\ntwo) \times 2$ is $2+\cdots+2+2$,
and is thus connected to $S(n-1)+1$ by a $1+1 \leftrightarrow 2$ move.

For $n \equiv 2 \bmod 4$, $S(n-1)+1$ is again $S(n-2)+1+1$, but 
now $n-2 \equiv 0 \bmod 4$, so this is $S(\ntwo -1)\times2 +1+1$.
On the other hand, $S(\ntwo)\times2$ is $S(\ntwo -1)\times2+2$ 
(because $\ntwo$ is odd), so $S(n-1)+1$ and $S(\ntwo)\times2$ are 
also connected by a $1+1 \leftrightarrow 2$ move, and we are done.
\qed\bigskip

As $n$ increases, the head of $\B(n)$ remains unchanged aside from
adding $+1$s to each partition.  So there is a single infinite
sequence $\B=\B_1,\B_2,\ldots$, beginning
$$
\emptyset, 2, 22, 4, 42, 222, 2222, 422, 44, 8,
82, 442, 4222, 22222, 222222, 42222, 4422, 822, 84, 444, \ldots
$$
such that each $\B(n)$ is just the initial substring of $\B$ of
partitions summing to $\leq n$, padded with the appropriate number 
of $1$s.  

$\B$ is a list of all binary partitions with only even parts, so
halving each one gives a Gray sequence $\B/2$ of the binary partitions
of all $n$.  Here the notion of a legal transition must be expanded to
include the operation $P \to P+1$, the remnant of the transition $1+1
\to 2$ after dropping all 1s and then halving.  By the construction,
the subsequence of partitions in $\B/2$ with constant sum $n$ is
identical to $\B(n)$ above if $n$ is even, and is the reverse of
$\B(n)$ if $n$ is odd.

\section{Stepping through the sequence}

Given a partition $P$ in the sequence $\B$, we can calculate the
partition which comes before or after $P$ easily.  We will give
explicit maps $\phi_+$ and $\phi_-$ which take a binary partition
and return which of the rules $2^k+2^k\leftrightarrow 2^{k+1}$
transforms $P$ into its successor or predecessor in $\B$.

Looking back at the construction of $\B$ from the $\B(n)$, we see that
for most $P$, the adjacent partitions $\phi_\pm(P)$ are the same size
as $P$; by ``size'' or $|P|$ we mean the sum of the parts of $P$ when
viewed as an element of $\B$, so after discarding any parts of size
$1$.  The only exception is when we try to apply $\phi_+$ to the last
partition of size $n$, which we called $S(n)$, or to apply $\phi_-$ to
the first partition in $\B$ of size $n$, which we now name $R(n)$.  In
the construction of $\B(n)$ this is the first partition which results
from a $\times2$ operation.

From the construction of the $\B(n)$, we can calculate:
\begin{equation}
\label{RS}
\begin{array}{rrcl}
\mbox{For } n\equiv  0\bmod 4, & R(n) &=& 2+\cdots+2, \\
\mbox{For } n\equiv  2\bmod 4, & S(n) &=& 2+\cdots+2, \\
\mbox{For } n\equiv  0\bmod 4, & S(n) &=& S(n/2) \times 2 = S(n/4) \times 4 = \ldots \\
&&=& \underbrace{2^a+\cdots+2^a}_b, \mbox{ where $n=2^a b$ with $b$ odd},\\
\mbox{For } n\equiv  2\bmod 4, & R(n) &=& S(n-2)+2.
\end{array}
\end{equation}
Conveniently, $S(n)$ for $n\equiv 2\bmod 4$ fits the $0\bmod 4$ pattern as well.

So the special cases are easily identified, and we already know
which transformation rules to apply to them:
\begin{equation}
\label{phi_base}
\begin{array}{rcl}
\phi_+(S(n)) &=&  1+1 \rightarrow 2 \\
\phi_-(R(n)) &=&  2 \rightarrow 1+1
\end{array}
\end{equation}
Of course, all parts of size $1$ are suppressed in $\B$, so the effect
of the rule $1+1\rightarrow 2$ is that a part of size $2$ appears from
nowhere, and the apparent size of the partition increases.

Aside from these special cases, all other transformations can be
determined recursively:
\begin{equation}
\label{phi_recursion}
\phi_{\pm}(P) = \left\{ \begin{array}{ll}
2\times \phi_\pm(\lfloor P/2 \rfloor), & |P| \equiv 0 \bmod 4 \\
2\times \phi_\mp(\lfloor P/2 \rfloor), & |P| \equiv 2 \bmod 4
\end{array}\right.
\end{equation}
Here $\lfloor \cdot \rfloor$ denotes deleting all parts of size $1$ to
again get an even binary partition, and the $\mp$ in the $2 \bmod 4$
case accounts for the reversal of $\B(\ntwo)$ in the definition of
$\B(n)$.  The multiplication by 2 acts on the rule returned, so
$2\times (1+1 \rightarrow 2)$ is $2+2 \rightarrow 4$.

Unravelling the recursion leads to a startlingly quick way to step
through $\B$.

\begin{thm}
\label{thm_phi}
Suppose you are given a binary partition $P$ written 
$\ldots d_3 d_2 d_1 d_0$, where the the digit $d_k$ is the number of
parts of $P$ of size $2^k$.  Since we are working in $\B$, we will
ignore the value of the $1$s place $d_0$, except to assume it is
at least two if the transformation $1+1\rightarrow2$ is needed.
\begin{itemize}
\item
Let $i$ be maximal with $d_i>0$  ($i=0$ if $P=\emptyset$).
\item
Let $j$ be the second largest integer with $d_j>0$ (or $j=0$ if none).
\item
Let $\epsilon$ be $(-1)^{\sum d_k}$ summing over $1\leq k\leq i-1$.
\end{itemize}
Then the transformations $\phi_\pm$ act on $P$ according to the
following rules:
\newcommand{\incr}{\!\uparrow}
\newcommand{\decr}{\!\downarrow}
\newcommand{\Incr}{\incr\incr}
\newcommand{\Decr}{\decr\decr}
$$
\newlength{\gnat}\settowidth{\gnat}{odd}
\newcommand{\case}[1]{\mbox{\rm ({#1})} & }
\begin{array}{lrl}
&\underline{\mbox{If $(d_i,d_j,\epsilon)$ is\ldots}}
&\underline{\mbox{then $\phi_\pm$ says\ldots}} \\[6pt]
\case{a}  
(\makebox[\gnat]{$1$},*,\mp1) & d_i \decr, d_{i-1} \Incr,
\mbox{ split a largest part}\\
\case{b}  \mbox{any other }
(\mbox{odd},*,\mp1) & d_{i+1} \incr, d_i \Decr,
\mbox{ merge two largest parts}\\[4pt]
\case{c} 
(\mbox{odd},1,\pm1) & d_j \decr, d_{j-1} \Incr,
\mbox{ split a $2^{nd}$-largest part}\\
\case{d}  \mbox{any other }
(\mbox{odd},*,\pm1) & d_{j+1} \incr, d_j \Decr,
\mbox{ merge two $2^{nd}$-largest parts}\\[4pt]
\case{e}
(\mbox{even},*,\mp1) & d_i \decr, d_{i-1} \Incr,
\mbox{ split a largest part}\\[4pt]
\case{f}
(\mbox{even},*,\pm1) & d_{i+1} \incr, d_i \Decr,
\mbox{ merge two largest parts}
\end{array}
$$
Here $\,\incr$, $\decr$ indicate an increment or decrement by
one, and $\,\Incr$, $\Decr$ by two.
\end{thm}
Note in particular that every transition in $\B$ involves merging or
splitting one of the largest two sizes of parts in $P$!  As a result,
successive updates can be done in place in constant time, in a
well-chosen data type (where you keep track of $\epsilon$, and 
you never need to search for $d_j$, e.g. a linked list of nonzero $d_k$.)

\proof
The recursive part of the definition of $\phi_\pm$ is handled trivially in this
notation.  The operation $\lfloor P \rfloor$ just shifts each digit to the right
and forgets $d_0$, and the value of $|P| \bmod 4$ depends only on the parity
of the 2s digit $d_1$.  So we keep dropping rightmost digits until we arrive 
at one of the base cases $\phi_+(S(n))$ or $\phi_-(R(n))$, and just remember 
$\epsilon$ to know if we've reversed direction an even or odd number of 
times\footnote{The fact that we ignore $d_0$ in the definition of $\epsilon$ 
perhaps reflects a moral imperfection in our definition of $\B(n)$: maybe for odd 
$n$ it should be the reverse of $\B(n-1)$, which then eliminates the reversal in 
the $2\bmod 4$ case.}.

The rest of the proof consists of identifying which of the base cases listed in 
equations~(\ref{RS}) and~(\ref{phi_base}) is the destination of
each partition:
\begin{itemize}
\item
Rule (d) is for partitions that end at $\phi_+(S(n))$.  The (odd) lead digit $d_i$ is
the (odd) number of identical parts in $S(n)=2^a+\cdots+2^a$.  We reach the base
case as soon as the second-largest nonzero digit $d_j$ is deleted, and the 
transformation $1+1\rightarrow2$ therefore joins two $2^j$s into a $2^{j+1}$.
\item
Rule (c) takes care of the exception where we run into $\phi_-(R(2n+2))$ one
step before we would otherwise reach $\phi_+(S(n))$.
\item
Rule (e) is for partitions that end at $\phi_-(R(n))$ when $n\equiv 0\bmod 4$,
an even-length sum $2+\cdots+2$.
\item
Rule (a) is for the somewhat special case $\phi_-(R(2))$.  $R(2)=2$, the only 
time $R(n)$ is a partition with an odd number of parts, since $S(0)=\emptyset$.
\item
Rules (b) and (f) are the fall-through cases: the recursion avoids all base cases
until it gets to $\phi_+(\emptyset)$, to which the transformation $1+1\rightarrow2$
is applied.
\end{itemize}
\qed

For example, let us compute the next several partitions in $\B$ beginning
with $P=256^5 32^2 16^1 4^4 2^3$ (so $|P|=1382$).
$$
\renewcommand{\ss}{\scriptscriptstyle}
\begin{array}{rrrrrrrrrrrcl}
\ss 1024 & \ss 512 & \ss 256 & \ss 128 & \ss 64 & \ss 32 & \ss 16 & \ss 8 & \ss 4 & \ss 2 & \ss 1 
& \quad\epsilon\quad \\
 & & 5 & 0 & 0 & 2 & 1 & 0 & 4 & 3 & 0 &   +1 & \mbox{Rule (d)} \\
 & & 5 & 0 & 1 & 0 & 1 & 0 & 4 & 3 & 0 &   -1 & \mbox{Rule (b)} \\
 & 1 & 3 & 0 & 1 & 0 & 1 & 0 & 4 & 3 & 0 &   +1 & \mbox{Rule (d)} \\
 & 2 & 1 & 0 & 1 & 0 & 1 & 0 & 4 & 3 & 0 &   +1 & \mbox{Rule (f)} \\
1 & 0 & 1 & 0 & 1 & 0 & 1 & 0 & 4 & 3 & 0 &   +1 & \mbox{Rule (c)} \\
1 & 0 & 0 & 2 & 1 & 0 & 1 & 0 & 4 & 3 & 0 &   -1 & \mbox{Rule (a)} \\
 & 2 & 0 & 2 & 1 & 0 & 1 & 0 & 4 & 3 & 0 &   -1 & \mbox{Rule (e)} \\
 & 1 & 2 & 2 & 1 & 0 & 1 & 0 & 4 & 3 & 0 &   -1 & \quad\ldots
\end{array}
$$
Note that $\epsilon$ changes sign under rules~(b) and~(c), and under (d) unless $i=j+1$.

\section{Calculating individual terms}

The recursive definition leading to the sequence $\B=\B_1,\B_2,\ldots$
allows us to explicitly compute the bijection $k \leftrightarrow \B_k$.  First
we introduce an alternate notation for even binary partitions.

\begin{defn}
Let $P$ be an even binary partition (that is, with no parts of size 1).
Define the {\em trail} of $P$, $\tau(P)=\tau=\tau_0,\tau_1,\tau_2,\ldots$ by
$$
\tau_i = \left | \lfloor P/2^i \rfloor \right |
$$
where again $\lfloor \cdot \rfloor$ indicates deleting any parts of size 1 
or smaller to again get an even binary partition.
\end{defn}

In other words, $\tau_i$ is the size of the partition after $i$
iterations of the map ``halve all parts and delete parts of size 1.''
For example, if $P=88422$, then $\lfloor P/2 \rfloor = 442$ and
$\lfloor P/4 \rfloor = 22$, so $\tau(P)=24,10,4,0,\ldots$, and we
henceforth omit the trailing 0s.  The partition $P$ is easily recovered
from $\tau(P)$: the number of parts of $P$ of size $2^i$ is 
$\tau_{i-1}/2 - \tau_i$, the number of parts of size 1 dropped on the $i$th 
iteration.

\begin{thm}
\label{thm_Bk}
Given just the integer $k$, the trail of $\B_k$ can be determined as follows:
\begin{enumerate}\renewcommand{\theenumi}{\roman{enumi}}
\item
$|\B_k|$ is the smallest $n$ such that $k \leq b(n)$.
\item
With $n=|\B_k|$ as above, $\lfloor \B_k/2 \rfloor = \B_\ell$, where
$$
\begin{array}{ll}
\ell = k - b(n-2), & n \equiv 0 \bmod 4, \\[4pt]
\ell = b(n) + 1 - k, & n \equiv 2 \bmod 4.
\end{array}
$$
\end{enumerate}

Conversely, the trail $\tau = \tau_0,\tau_1,\tau_2,\ldots$ corresponds
to the $k$th partition $\B_k$ if its truncation $\tau_1,\tau_2,\tau_3,\ldots$
corresponds to the $\ell$th, where
$$
\begin{array}{ll}
k = b(\tau_0-2) + \ell & \tau_0 \equiv 0 \bmod 4, \\[4pt]
k = b(\tau_0) + 1 - \ell, & \tau_0 \equiv 2 \bmod 4.
\end{array}
$$
This recurrence and the base case $\tau=0,0,\ldots$ corresponding to
$\B_{b(0)} = \B_1 =\emptyset$ suffices to determine the location in $\B$ of the
binary partition with any given trail.
\end{thm}

For example, partition $88422$ with trail $24,10,4$ appears at position 
$b(22) + b(10) + 1 - (b(2) + b(0)) = 86$.  And to find $\B_{123456789}$,
we calculate that
$$
\begin{array}{rcll}
123456789 &\in& (b(646),b(648)] & \mbox{note } 648\equiv 0\bmod 4, \\
123456789-b(646) &\in& (b(304),b(306)] & \mbox{note } 306\equiv 2\bmod 4, \\
b(306)+1-(123456789-b(646)) &\in& (b(120),b(122)] & \mbox{etc,}
\end{array}
$$
and its trail is $648, 306, 122, 58, 28, 14$, so the partition
is $64^7 32^0 16^1 8^3 4^{31} 2^{18}$.  This amounts to a writing of 123456789 
in terms of values of $b(n)$ with $n\equiv 2\bmod 4$:
$$
b(646) + b(306) + 1 - \left( b(122) + 1 -\left( b(58) + 1 -
  \left(b(26)+b(14)+1-\left(b(0) \rule{0pt}{10pt}\right)\!\right)\!\right)\!\right)
$$
It is possible to work with this representation directly, using sign changes to
track the ${}\bmod 4$ behavior, but statements tend to be inelegant.

Both maps rely heavily on the values of the function $b(n)$ which counts
the number of binary partitions of $n$, and for which we have no closed form.
From a computational complexity point of view this is inevitable, since from
the map $k\to\B_k$ the values of $b(n)$ can be determined easily.  In practice
the computation of $b(n)$ by recurrence~(\ref{bn}) is inexpensive.

\proof
Since the trail of any even binary partition $P$ begins with $\tau_0=|P|$,
the given maps $k\to\tau$ and $\tau\to k$ are clearly inverse, so it suffices to
show either direction.

The map $k\to\tau(\B_k)$ follows directly from the construction of the 
sequences $\B(n)$.  Suppose $n$ is even.  After deleting all parts of size 1
from the partitions in all $\B(n)$, the first $b(n-2)$ terms of $\B(n)$ are 
exactly $\B(n-2)$, and the remaining $b(n/2)$ terms are the even binary 
partitions of $n$.  So these partitions occur in $\B$ as a solid block of terms 
with indices in $(b(n-2),b(n)]$, justifying~(i).  

For~(ii), recall that each of the terms in this block was obtained by doubling 
some term in $\B(n/2)$; the $\lfloor\cdot\rfloor$ operation deletes all parts of 
size 1 and allows us to work in $\B$ instead.  If $n\equiv 0\bmod 4$ then we
doubled $\B_\ell$ to get the $\ell$th term in the block, $\B_{b(n-2)+\ell}$, while
if $n\equiv 2\bmod 4$ then $\B(n/2)$ was reversed, and we doubled $\B_\ell$ 
to get the $\ell$th term from the end, $\B_{b(n)+1-\ell}$.
\qed

\end{document}